\documentclass{article}[12pt]

\usepackage{latexsym}
\usepackage{amsmath}
\usepackage{amssymb}

\begin{document}

\title{Explicitly solvable algebraic equations of degree 8 and 9}

\author{Francesco Calogero$^{a,b}$\thanks{e-mail: francesco.calogero@roma1.infn.it}
\thanks{e-mail: francesco.calogero@uniroma1.it}
 , Farrin Payandeh$^c$\thanks{e-mail: farrinpayandeh@yahoo.com}
 \thanks{e-mail: f$\_$payandeh@pnu.ac.ir}}

\maketitle   \centerline{\it $^{a}$Physics Department, University of
Rome "La Sapienza", Rome, Italy}

\maketitle   \centerline{\it $^{b}$INFN, Sezione di Roma 1}

\maketitle

\maketitle   \centerline{\it $^{c}$Department of Physics, Payame
Noor University, PO BOX 19395-3697 Tehran, Iran}

\maketitle

\begin{abstract}

The \textit{generic} monic polynomial of degree $N$ features $N$ \textit{a
priori arbitrary} coefficients $c_{m}$ and $N$ zeros $z_{n}$. In this paper
we limit consideration to $N=8$ and $N=9$. We show that if the $N$---\textit{%
a priori arbitrary}---coefficients $c_{m}$ of these polynomials are \textit{%
appropriately }defined---as it were, \textit{a posteriori}---in terms of $6$
\textit{arbitrary} parameters, then the $N$ roots of these polynomials can
be \textit{explicitly} computed in terms of radicals of these $6$
parameters. We also report the \textit{constraints }on the $N$ coefficients $%
c_{m}$ implied by the fact that they are so defined in terms of $6$ \textit{%
arbitrary} parameters; as well as the \textit{explicit} determination of
these $6$ parameters in terms of the $N$ coefficients $c_{m}$.

\end{abstract}

\section{Introduction}

This paper is a follow-up to the paper \cite{CP2021}, and we refer to that
paper for the motivation of this kind of investigations.

\textbf{Notation 1-1}. Hereafter the $N$ zeros $z_{n}$ and the $N$
coefficients $c_{m}$ of the monic polynomial $P_{N}\left( z\right) $ of
degree $N$ are defined as follows:
\begin{equation}
P_{N}\left( z\right) =\prod\limits_{n=1}^{N}\left( z-z_{n}\right)
=z^{N}+\sum_{m=0}^{N-1}\left( c_{m}z^{m}\right) ~;  \label{PNz}
\end{equation}%
and in this paper we limit consideration to $N=8$ and $N=9$. $\blacksquare $

In this paper we identify \textit{special} cases of the polynomial (\ref{PNz}%
) with $N=8$ and $N=9$---with its $N$ coefficients $c_{m}$ \textit{%
appropriately} defined in terms of $6$ \textit{arbitrary} parameters---so
that its $N$ zeros $z_{n}$ can then be themselves \textit{explicitly}
determined in terms of \textit{radicals} of these $6$ parameters; and we
also report below both the expressions of the $6$ parameters in terms of the
$N$ coefficients $c_{m}$ and the \textit{constraints} on the coefficients $%
c_{m}$ implied by this approach.

Of course the results reported below do not contradict the implications of
the findings about the zeros of polynomials implied by the Theory of Galois
Groups, and are certainly contained as special cases within that framework;
but they are obtained by much more elementary means.

In the following two Sections we report separately our findings for the two
cases with $N=8$ and $N=9$; the proofs of these findings are provided in
\textbf{Appendix A}. A terse final Section outlines possible future
developments.

\bigskip

\section{Results for N=8}

\textbf{Proposition 2-1}. Assume that the $8$ coefficients $c_{m}$ of the
polynomial (\ref{PNz}) with $N=8$ may be expressed as follows in terms of
the $6$ \textit{arbitrary} parameters $\alpha _{0},$ $\alpha _{1},$ $\beta
_{0},$ $\beta _{1},$ $\gamma _{0},$ $\gamma _{1}$:
\begin{subequations}
\label{cc}
\begin{equation}
c_{7}=4\alpha _{1}~,  \label{c7}
\end{equation}%
\begin{equation}
c_{6}=4\alpha _{0}+6\left( \alpha _{1}\right) ^{2}+2\beta _{1}~,  \label{c6}
\end{equation}%
\begin{equation}
c_{5}=12\alpha _{0}\alpha _{1}+4\left( \alpha _{1}\right) ^{3}+6\alpha
_{1}\beta _{1}~,  \label{c5}
\end{equation}%
\begin{equation}
c_{4}=6\alpha _{0}\left[ \alpha _{0}+2\left( \alpha _{1}\right) ^{2}+\beta
_{1}\right] +\left( \alpha _{1}\right) ^{2}\left[ \left( \alpha _{1}\right)
^{2}+6\beta _{1}\right] +2\beta _{0}+\left( \beta _{1}\right) ^{2}+\gamma
_{1}~,  \label{c4}
\end{equation}%
\begin{equation}
c_{3}=4\alpha _{0}\alpha _{1}\left[ 3\alpha _{0}+\left( \alpha _{1}\right)
^{2}+3\beta _{1}\right] +2\alpha _{1}\left[ 2\beta _{0}+\left( \alpha
_{1}\right) ^{2}\beta _{1}+\left( \beta _{1}\right) ^{2}+\gamma _{1}\right]
~,  \label{c3}
\end{equation}%
\begin{eqnarray}
c_{2} &=&2\alpha _{0}\left\{ 2\left( \alpha _{0}\right) ^{2}+3\alpha
_{0}\left( \alpha _{1}\right) ^{2}+2\beta _{0}+3\left[ \alpha _{0}+\left(
\alpha _{1}\right) ^{2}\right] \beta _{1}+\left( \beta _{1}\right)
^{2}\right\}  \nonumber \\
&&+2\left( \alpha _{1}\right) ^{2}\beta _{0}+\left[ 2\beta _{0}+\left(
\alpha _{1}\right) ^{2}\beta _{1}\right] \beta _{1}+\left[ 2\alpha
_{0}+\left( \alpha _{1}\right) ^{2}+\beta _{1}\right] \gamma _{1}~,
\label{c2}
\end{eqnarray}%
\begin{equation}
c_{1}=\alpha _{1}\left\{ 2\alpha _{0}\left[ 2\left( \alpha _{0}\right)
^{2}+2\beta _{0}+\left( 3\alpha _{0}+\beta _{1}\right) \beta _{1}\right]
+2\beta _{0}\beta _{1}+\left( 2\alpha _{0}+\beta _{1}\right) \gamma
_{1}\right\} ~,  \label{c1}
\end{equation}%
\begin{eqnarray}
c_{0} &=&\alpha _{0}\left\{ \left( \alpha _{0}\right) ^{3}+\alpha _{0}\left[
2\beta _{0}+2\left( \alpha _{0}\right) ^{2}\beta _{1}+\left( \beta
_{1}\right) ^{2}\right] +2\beta _{0}\beta _{1}+\left( \alpha _{0}+\beta
_{1}\right) \gamma _{1}\right\}  \nonumber \\
&&+\beta _{0}\left( \beta _{0}+\gamma _{1}\right) +\gamma _{0}~.  \label{cc0}
\end{eqnarray}

Then the $8$ roots $z_{n}\equiv z_{\lambda \mu \nu }$ ($n=1,2,3,4,5,6,7,8$; $%
\lambda =0,1,~\mu =0,1,~\nu =0,1$) of the polynomial (\ref{PNz}) with $N=8$
are \textit{explicitly} given, in terms of the $6$ parameters $\alpha _{0},$
$\alpha _{1},$ $\beta _{0},$ $\beta _{1},$ $\gamma _{0},$ $\gamma _{1}$, by
the following definitions: the $8$ numbers $z_{\lambda \mu \nu }$ are the $2$
roots (with $\lambda =1,2,$ $\mu =1,2$, $\nu =1,2$) of the following $4$
\textit{quadratic} equation in $z$,
\end{subequations}
\begin{subequations}
\label{Eqzyx}
\begin{equation}
z^{2}+\alpha _{1}z+\alpha _{0}=y_{\mu \nu }~,~~~\mu =1,2,~~~\nu =1,2~,
\label{Equaz}
\end{equation}%
where the $4$ numbers $y_{\mu \nu }$ (with $\mu =1,2$, $\nu =1,2$) are the $%
2 $ roots of the following $2$ \textit{quadratic }equations in $y$,
\begin{equation}
y^{2}+\beta _{1}y+\beta _{0}=x_{\nu }~,~~~\nu =1,2~,  \label{Equay}
\end{equation}%
where the $2$ numbers $x_{\nu }$ (with $\nu =1,2$) are the $2$ roots of the
following \textit{quadratic }equation in $x$,%
\begin{equation}
x^{2}+\gamma _{1}x+\gamma _{0}=0~.~~~\blacksquare  \label{Equax}
\end{equation}

\textbf{Remark 2-1}. The $3$ \textit{quadratic} eqs. (\ref{Eqzyx}) can of
course be solved \textit{explicitly}. The resulting \textit{explicit}
formula expressing the $8$ zeros $z_{n}\equiv z_{\lambda \mu \nu }$ in terms
of the $6$ parameters $a_{0},$ $a_{1},$ $a_{2},$ $b_{0},$ $b_{1}$
involves---in a nested way---only \textit{square} roots; anybody interested
can obtain it easily, the only tool needed to that end is the formula giving
the $2$ roots $x_{\pm }$ of a monic polynomial of \textit{second} degree:
\end{subequations}
\begin{equation}
x^{2}+a_{1}x+a_{0}=\left( x-x_{+}\right) \left( x-x\right) ~,~~~x_{\pm
}=\left( -a_{1}\pm \sqrt{\left( a_{1}\right) ^{2}-4a_{0}}\right)
/2~.~~~\blacksquare
\end{equation}

Next, let us face the "inverse" task to assign---as it were, \textit{a priori%
}---the $8$ \textit{coefficients} $c_{m}$ of a monic polynomial of degree $8$
(see (\ref{PNz}) with $N=8$) and to then find the $6$ \textit{parameters} $%
\alpha _{0},$ $\alpha _{1},$ $\beta _{0},$ $\beta _{1},$ $\gamma _{0},$ $%
\gamma _{1}$ which determine---as it were, \textit{a posteriori}---via
\textit{explicit} formulas both these $8$ coefficients $c_{m}$ and the $8$
\textit{zeros} $z_{n}$ of that polynomial (\ref{PNz}), as well as \textit{%
explicit} formulas displaying the corresponding \textit{constraints} implied
by these assignments on the $8$ coefficients $c_{m}$.

The following proposition provides these findings, which complement those
reported in \textbf{Proposition 2-1}.

\textbf{Proposition 2-2}. If the $8$ parameters $c_{m}$ are expressed in
terms of the $6$ parameters $\alpha _{0},$ $\alpha _{1},$ $\beta _{0},$ $%
\beta _{1},$ $\gamma _{0},$ $\gamma _{1}$ by the $8$ formulas (\ref{cc}),
then the $4$ parameters $\alpha _{1},$ $\beta _{1},$ $\gamma _{0},$ $\gamma
_{1}$ are themselves expressed as follows in terms of the $8$ coefficients $%
c_{m}$ and of the "free" parameters $\alpha _{0}$ and $\beta _{0}$:
\begin{subequations}
\label{alfbet}
\begin{equation}
\alpha _{1}=c_{7}/4~,  \label{alf1}
\end{equation}%
\begin{equation}
\beta _{1}=-\left[ 32\alpha _{0}-8c_{6}+3\left( c_{7}\right) ^{2}\right]
/16~,  \label{bet1}
\end{equation}%
\begin{eqnarray}
&\gamma _{0}=c_{0}-\left\{ \alpha _{0}\left[ 8c_{6}-3\left( c_{7}\right) ^{2}%
\right] -16\left[ \left( \alpha _{0}\right) ^{2}-\beta _{0}\right] \right\}
\cdot &  \notag \\
&\cdot \left\{ \left( c_{7}\right) ^{4}-8\left( c_{6}\right)
^{2}+32c_{4}-2\alpha _{0}\left[ 8c_{6}-3\left( c_{7}\right) ^{2}\right] +32%
\left[ \left( \alpha _{0}\right) ^{2}-\beta _{0}\right] \right\} /512~,&
\label{gam0}
\end{eqnarray}%
\begin{equation}
\gamma _{1}=-\alpha _{0}c_{6}+2\left[ \left( \alpha _{0}\right) ^{2}-\beta
_{0}\right] +c_{4}+\left[ -8\left( c_{6}\right) ^{2}+12\alpha _{0}\left(
c_{7}\right) ^{2}+\left( c_{7}\right) ^{4}\right] /32~,  \label{gam1}
\end{equation}%
while the $8$ coefficients $c_{m}$ satisfy the following $4$ \textit{%
constraints}:
\end{subequations}
\begin{subequations}
\label{Con}
\begin{equation}
c_{5}=c_{7}\left[ 24c_{6}-7\left( c_{7}\right) ^{2}\right] /32~,
\label{Con1}
\end{equation}%
\begin{equation}
c_{3}=c_{7}\left[ 128c_{4}-20c_{6}\left( c_{7}\right) ^{2}+7\left(
c_{7}\right) ^{4}\right] /256~,  \label{Con2}
\end{equation}%
\begin{eqnarray}
c_{2} &=&\left\{ 512c_{4}\left[ 4c_{6}-\left( c_{7}\right) ^{2}\right] -64%
\left[ 8c_{6}-3\left( c_{7}\right) ^{2}\right] \left( c_{6}\right)
^{2}\right.  \nonumber \\
&&\left. +\left[ 16c_{6}-7\left( c_{7}\right) ^{2}\right] \left(
c_{7}\right) ^{4}\right\} /4096~,  \label{Con3}
\end{eqnarray}%
\begin{equation}
c_{1}=c_{7}\left[ 8c_{6}-3\left( c_{7}\right) ^{2}\right] \left[
32c_{4}-8\left( c_{6}\right) ^{2}+\left( c_{7}\right) ^{4}\right]
/2048~.~~~\blacksquare  \label{Con4}
\end{equation}

\bigskip

\section{Results for N=9}

The findings reported in this \textbf{Section 3} are \textit{analogous}, but
\textit{different}, from those reported in \textbf{Section 2}; accordingly,
the variables and parameters used in this \textbf{Section 3 }are \textit{%
different} from those having the \textit{same} names in \textbf{Section 2},
although they play \textit{analogous} roles.

\textbf{Proposition 3-1}. Assume that the $9$ coefficients $c_{m}$ of the
polynomial (\ref{PNz}) with $N=9$ may be expressed as follows in terms of
the $6$ \textit{arbitrary} parameters $\alpha _{0},$ $\alpha _{1},$ $\alpha
_{2},$ $\beta _{0},$ $\beta _{1},$ $\beta _{2}$:
\end{subequations}
\begin{subequations}
\label{31cc}
\begin{equation}
c_{0}=\alpha _{0}\left[ \beta _{1}+\alpha _{0}\left( \alpha _{0}+\beta
_{2}\right) \right] +\beta _{0}~,  \label{31c0}
\end{equation}%
\begin{equation}
c_{1}=\alpha _{0}\alpha _{1}\left( 3\alpha _{0}+2\beta _{2}\right) +\alpha
_{1}\beta _{1}~,  \label{31c1}
\end{equation}%
\begin{equation}
c_{2}=3\alpha _{0}\left[ \left( \alpha _{1}\right) ^{2}+\alpha _{0}\alpha
_{2}\right] +\beta _{2}\left[ \left( \alpha _{1}\right) ^{2}+2\alpha
_{0}\alpha _{2}\right] +\alpha _{2}\beta _{1}~,  \label{31c2}
\end{equation}%
\begin{equation}
c_{3}=3\alpha _{0}\left( \alpha _{0}+2\alpha _{1}\alpha _{2}\right) +2\beta
_{2}\left( \alpha _{0}+\alpha _{1}\alpha _{2}\right) +\left( \alpha
_{1}\right) ^{3}+\beta _{1}~,  \label{31c3}
\end{equation}%
\begin{equation}
c_{4}=3\alpha _{0}\left[ 2\alpha _{1}+\left( \alpha _{2}\right) ^{2}\right]
+3\left( \alpha _{1}\right) ^{2}\alpha _{2}+\left[ 2\alpha _{1}+\left(
\alpha _{2}\right) ^{2}\right] \beta _{2}~,  \label{31c4}
\end{equation}%
\begin{equation}
c_{5}=3\alpha _{1}\left[ \alpha _{1}+\left( \alpha _{2}\right) ^{2}\right]
+2\alpha _{2}\left( 3\alpha _{0}+\beta _{2}\right) ~,  \label{31c5}
\end{equation}%
\begin{equation}
c_{6}=3\alpha _{0}+\alpha _{2}\left[ 6\alpha _{1}+\left( \alpha _{2}\right)
^{2}\right] +\beta _{2}~,  \label{31c6}
\end{equation}%
\begin{equation}
c_{7}=3\left[ \alpha _{1}+\left( \alpha _{2}\right) ^{2}\right] ~,
\label{31c7}
\end{equation}%
\begin{equation}
c_{8}=3\alpha _{2}~.  \label{31c8}
\end{equation}

Then the $9$ roots $z_{n}\equiv z_{\lambda \mu }$ (with $n=1,..,9$ and $%
\lambda =1,2,3$, $\mu =1,2,3$) of the polynomial (\ref{PNz}) with $N=9$ are
\textit{explicitly} given, in terms of the $6$ parameters $\alpha _{0},$ $%
\alpha _{1},$ $\alpha _{2},$ $\beta _{0},$ $\beta _{1},$ $\beta _{2}$ by the
following definitions: $z_{n}\equiv z_{\lambda \mu }$ is one of the $3$
roots of the following $3$ \textit{cubic }equations in $z$,
\end{subequations}
\begin{subequations}
\label{Eqyz}
\begin{equation}
z^{3}+\alpha _{2}z^{2}+\alpha _{1}z+\alpha _{0}=y_{\mu }~,~~~\mu =1,2,3~,
\label{31Eqz}
\end{equation}%
where the $3$ numbers $y_{\mu }\ $are the $3$ roots of the following \textit{%
cubic} equation in $y$:
\begin{equation}
y^{3}+\beta _{2}y^{2}+\beta _{1}y+\beta _{0}=0~.~\blacksquare  \label{31Eqy}
\end{equation}

\textbf{Proposition 3-2}. If the $9$ parameters $c_{m}$ of the polynomial (%
\ref{PNz}) with $N=9$ are expressed in terms of the $6$ parameters $\alpha
_{0},$ $\alpha _{1},$ $\alpha _{2},$ $\beta _{0},$ $\beta _{1},$ $\beta _{2}$
by the $9$ formulas (\ref{31cc}), then the $6$ parameters $\alpha _{0},$ $%
\alpha _{1},$ $\alpha _{2},$ $\beta _{0},$ $\beta _{1},$ $\beta _{2}$ are
themselves expressed as follows in terms of the $9$ coefficients $c_{m}$:
\end{subequations}
\begin{subequations}
\label{32alfbet}
\begin{equation}
\alpha _{2}=c_{8}/3~,  \label{32alf2}
\end{equation}%
\begin{equation}
\alpha _{1}=\left[ 3c_{7}-\left( c_{8}\right) ^{2}\right] /9~,
\label{32alf1}
\end{equation}%
\begin{equation}
\alpha _{0}=\frac{c_{4}-3\left( \alpha _{1}\right) ^{2}\alpha _{2}}{3\left[
2\alpha _{1}+\left( \alpha _{2}\right) ^{2}\right] }-\beta _{2}/3~,
\label{32alf0}
\end{equation}%
\begin{equation}
\beta _{2}=\frac{c_{3}-3\alpha _{0}\left( \alpha _{0}-2\alpha _{1}\alpha
_{2}\right) -\left( \alpha _{1}\right) ^{3}-\beta _{1}}{2\left( \alpha
_{0}+\alpha _{1}\alpha _{2}\right) }~,  \label{32bet2}
\end{equation}%
\begin{equation}
\beta _{1}=\left[ c_{1}-\alpha _{0}\alpha _{1}\left( 3\alpha _{0}+2\beta
_{2}\right) \right] /\alpha _{1}~,  \label{32bet1}
\end{equation}%
\begin{equation}
\beta _{0}=c_{0}-\alpha _{0}\left[ \left( \alpha _{0}\right) ^{2}+\alpha
_{0}\beta _{2}+\beta _{1}\right] ~;  \label{32bet0}
\end{equation}%
note that these formulas (\ref{32alfbet}), via their sequentially nested
character, express \textit{explicitly} the $6$ quantities $\alpha _{0},$ $%
\alpha _{1},$ $\alpha _{2},$ $\beta _{0},$ $\beta _{1},$ $\beta _{2}$ in
terms of the $6$ coefficients $c_{8},$ $c_{7},$ $c_{4},$ $c_{3},$ $c_{1},$ $%
c_{0}$.

While the $9$ coefficients $c_{m}$ are themselves required to satisfy as
many as $5$ \textit{constraints}, which can be expressed in the form of
\textit{explicit expressions} of any $5$ of the $6$ coefficients $c_{\ell }$
with $\ell =1,2,...,6$ in terms of the remaining coefficients $c_{m}$ with $%
m=1,2,...,8$: as explained in \textbf{Appendix A}. $\blacksquare $

\bigskip

\section{Outlook}

It is possible---but not necessarily useful---to try and extend the
simple-minded approach utilized in this paper in order to identify other
\textit{special} polynomials---of higher order than $9$, and characterized
by \textit{appropriate} restrictions on their coefficients---which allow
their zeros to be \textit{explicitly} computed in terms of their
coefficients; or even to treat in analogous manner polynomial equations the
arguments and coefficients of which are more complicated objects than
\textit{scalar} quantities (for instance, \textit{matrices}).

\bigskip

\section{Appendix A}

In this Appendix we justify the results reported in \textbf{Sections 2 }and
\textbf{3}.

To derive the results reported in \textbf{Propositions 2-1} one should
firstly replace $x$ in the \textit{quadratic} eq. (\ref{Equax}) with the
expression $y^{2}+\beta _{1}y+\beta _{0}$ (see eq. (\ref{Equay})) and
secondly replace, in the resulting \textit{quartic} equation for $y$, this
variable with the expression $z^{2}+\alpha _{1}z+\alpha _{0}$ (see eq. (\ref%
{Equaz})), obtaining thereby an equation of degree $8$ for the variable $z;$
then identify the $8$ coefficients of this equation---after having expanded
it in powers of $z$---with the $8$ coefficients $c_{m}$ of the polynomial $%
P_{8}\left( z\right) ,$ see eq. (\ref{PNz}). Likewise, to obtain the results
reported in \textbf{Proposition 3-1}, replace\textbf{\ }$y$ in the \textit{%
cubic} eq. (\ref{31Eqy}) with $z^{3}+\alpha _{2}z^{2}+\alpha _{1}z+\alpha
_{0}$ (see eq. (\ref{31Eqz})), expand the outcome in powers of $z$ and
identify the $9$ coefficients of the resulting polynomial of degree $9$ with
the coefficients $c_{m}$ of the polynomial $P_{9}\left( z\right) ,$ see (\ref%
{PNz}).

The results reported in \textbf{Propositions 2-2} are consequences of the
inversion of the formulas (\ref{alfbet}), performed in order to obtain the
parameters (identified by Greek letters) they feature in terms of the $8$
coefficients $c_{m}$; they have been obtained via a judicious choice of the
order in which these equations are solved, and they have then been checked
via the algebraic manipulation package \textbf{Mathematica}.

The procedure to obtain the results reported in \textbf{Proposition 3-2} is
a bit more tricky. The issue is of course to invert the algebraic eqs. (\ref%
{32alfbet}). It is convenient to consider them in the opposite order in
which they are listed. Then from the last $2$ (i. e. eqs. (\ref{31c8}) and (%
\ref{31c7})) we immediately get
\end{subequations}
\begin{equation}
\alpha _{2}=c_{8}/3~,~~~\alpha _{1}=\left[ 3c_{7}-\left( c_{8}\right) ^{2}%
\right] /9~.  \label{Aalf12}
\end{equation}%
Hence hereafter we shall use---for simplicity---the variables $\alpha _{1}$
and $\alpha _{2}$ in place of the variables $c_{7}$ and $c_{8}$.

It is then convenient to rewrite the $3$ eqs. (\ref{31c6}), (\ref{31c5}) and
(\ref{31c4}) as follows:
\begin{subequations}
\label{Acc}
\begin{equation}
3\alpha _{0}+\beta _{2}=c_{6}-\alpha _{2}\left[ 6\alpha _{1}+\left( \alpha
_{2}\right) ^{2}\right] ~,  \label{Ac6}
\end{equation}%
\begin{equation}
3\alpha _{0}+\beta _{2}=\left\{ c_{5}-3\alpha _{1}\left[ \alpha _{1}+\left(
\alpha _{2}\right) ^{2}\right] \right\} /\left( 2\alpha _{2}\right) ~,
\label{Ac5}
\end{equation}%
\begin{equation}
3\alpha _{0}+\beta _{2}=\left[ c_{4}-3\left( \alpha _{1}\right) ^{2}\alpha
_{2}\right] /\left[ 2\alpha _{1}+\left( \alpha _{2}\right) ^{2}\right] ~;
\label{Ac4}
\end{equation}%
and it is also easily seen that the $3$ eqs. (\ref{31c3}), (\ref{31c2}) and (%
\ref{31c1}) can also be combined to yield the following $3$ relations:%
\begin{equation}
3\alpha _{0}+\beta _{2}=\left( c_{2}\alpha _{1}-c_{1}\alpha _{2}\right)
/\left( \alpha _{1}\right) ^{3}~,  \label{Ac21}
\end{equation}%
\begin{equation}
3\alpha _{0}+\beta _{2}=\left[ c_{3}\alpha _{1}-c_{1}-\left( \alpha
_{1}\right) ^{4}\right] /\left[ 2\left( \alpha _{1}\right) ^{2}\alpha _{2}%
\right] ~,  \label{Ac31}
\end{equation}%
\begin{equation}
3\alpha _{0}+\beta _{2}=\left[ c_{3}\alpha _{2}-c_{2}-\left( \alpha
_{1}\right) ^{3}\alpha _{2}\right] /\left\{ \alpha _{1}\left[ 2\left( \alpha
_{2}\right) ^{2}-\alpha _{1}\right] \right\} ~.  \label{Ac32}
\end{equation}%
Since the left-hand sides of these $6$ equations (\ref{Acc}) are \textit{all
equal}, and the right-hand sides involve---via (\ref{Aalf12})---only the $8$
parameters $c_{m}$ (with $m=1,2,...,8$), this implies that these $8$
coefficients $c_{m}$ must satisfy, among themselves, $5$ \textit{constraints}%
; and it is also clear that it is quite easy to reformulate these $5$
\textit{constraints} in the form of $5$ \textit{explicit expressions} of any
$5$ of the $6$ coefficients $c_{\ell }$ with $\ell =1,2,...,6$ in terms of
the remaining coefficients $c_{m}$ (with $m=1,2,...,8$)---by just solving
the system of $5$ \textit{algebraic linear equations }for these coefficients
implied by the $6$ eqs. (\ref{Acc}).

\bigskip

\textbf{Acknowledgements}. We like to acknowledge with thanks $3$ grants,
facilitating our collaboration---mainly developed via e-mail exchanges---by
making it possible for FP to visit $3$ times (one of which in the future)
the Department of Physics of the University of Rome "La Sapienza": $2$
granted by that University, and one granted jointly by the Istituto
Nazionale di Alta Matematica (INdAM) of that University and by the
International Institute of Theoretical Physics (ICTP) in Trieste in the
framework of the ICTP-INdAM "Research in Pairs" Programme.\textbf{\ }%
We also\ like to thank Fernanda Lupinacci who, in these difficult
times---with extreme efficiency and kindness---facilitated all the
arrangements necessary for the presence of FP with her family in Rome.

\end{subequations}


\begin{thebibliography}{9}
\bibitem{CP2021} F. Calogero and F. Payandeh, "Two classes of explicitly
solvable sextic equations", arXiv:2104.03072 [math.DS] April 8, 2021.
\end{thebibliography}
\end{document}